# Euler sums and a prime numbers conundrum

Donal F. Connon

dconnon@btopenworld.com

12 March 2008

## Abstract

This note highlights an interesting connection between Euler sums of even weight and prime numbers.

In the enlightening work of Flajolet and Salvy [10] we find a treasure trove of both linear and non-linear Euler sums defined by

$$S_{\pi,q} = \sum_{n=1}^{\infty} \frac{H_n^{(\pi_1)} H_n^{(\pi_2)} \ldots H_n^{(\pi_k)}}{n^q}$$

where $\pi = (\pi_1, \ldots, \pi_k)$ is a partition of an integer $p$ into $k$ summands so that

$$p = \pi_1 + \ldots + \pi_k$$

and $\pi_1 \le \pi_1 \le \ldots \le \pi_k$ **.** The quantity $w = q + \pi_1 + \ldots + \pi_k$ is called the weight and $k$ is referred to as the degree.

A selection of the Euler sums with **even** weights reported in [10] is set out below:

(1) $$4\sum_{n=1}^{\infty} \frac{H_n^{(1)}}{n^3} = 5\varsigma(4)$$

(2) $$4\sum_{n=1}^{\infty} \frac{\left[H_n^{(1)}\right]^2}{n^2} = 17\varsigma(4)$$

(3) $$24\sum_{n=1}^{\infty} \frac{\left[H_n^{(1)}\right]^2}{n^4} = 97\varsigma(6) - 48\varsigma^2(3)$$

(4) $$4\sum_{n=1}^{\infty} \frac{H_n^{(3)}}{n^3} = 2\varsigma(6) + 2\varsigma^2(3)$$

(5) $$6\sum_{n=1}^{\infty}\frac{H_n^{(2)}}{n^4}=-2\varsigma(6)+6\varsigma^2(3)$$

(6) $$12\sum_{n=1}^{\infty}\frac{\left[H_n^{(2)}\right]^2}{n^2}=-43\varsigma(6)+24\varsigma^2(3)$$

Flajolet and Salvy [10, p.23] also reported that

(7) $$140\sum_{n=1}^{\infty}\frac{H_n^{(2)}}{n^8}+200\sum_{n=1}^{\infty}\frac{H_n^{(3)}}{n^7}=-330\varsigma(10)+280\varsigma(3)\varsigma(7)+160\varsigma^2(5)$$

(8) $$70\sum_{n=1}^{\infty}\frac{H_n^{(2)}}{n^8}-20\sum_{n=1}^{\infty}\frac{H_n^{(4)}}{n^6}=-227\varsigma(10)+140\varsigma(3)\varsigma(7)+100\varsigma^2(5)$$

and adding the above two equations results in

(9) $$210\sum_{n=1}^{\infty}\frac{H_n^{(2)}}{n^8}+200\sum_{n=1}^{\infty}\frac{H_n^{(3)}}{n^7}-20\sum_{n=1}^{\infty}\frac{H_n^{(4)}}{n^6}=-557\varsigma(10)+420\varsigma(3)\varsigma(7)+260\varsigma^2(5)$$

Subtraction results in

(10) $$70\sum_{n=1}^{\infty}\frac{H_n^{(2)}}{n^8}+200\sum_{n=1}^{\infty}\frac{H_n^{(3)}}{n^7}+20\sum_{n=1}^{\infty}\frac{H_n^{(4)}}{n^6}=-103\varsigma(10)+140\varsigma(3)\varsigma(7)+60\varsigma^2(5)$$

In 2002, Choi and Srivastava [3] proved that

(11)
$$16\sum_{n=1}^{\infty}\frac{\left[H_n^{(1)}\right]^4}{n^2}+32\sum_{n=1}^{\infty}\frac{H_n^{(1)}H_n^{(3)}}{n^2}-48\sum_{n=1}^{\infty}\frac{\left[H_n^{(1)}\right]^2H_n^{(2)}}{n^2}+96\sum_{n=1}^{\infty}\frac{H_n^{(1)}H_n^{(2)}}{n^3}=373\varsigma(6)+150\varsigma^2(3)$$

(12) $$48\sum_{n=1}^{\infty}\frac{H_n^{(1)}H_n^{(2)}}{n^3}=-101\varsigma(6)+120\varsigma^2(3)$$

(13) $$288\sum_{n=1}^{\infty}\frac{H_n^{(1)}H_n^{(2)}}{n^4}+64\sum_{n=1}^{\infty}\frac{H_n^{(1)}H_n^{(3)}}{n^3}-164\sum_{n=1}^{\infty}\frac{\left[H_n^{(1)}\right]^2H_n^{(2)}}{n^3}+24\sum_{n=1}^{\infty}\frac{\left[H_n^{(2)}\right]^2}{n^3}=$$

$$-2357\varsigma(7)+1972\varsigma(3)\varsigma(4)-88\varsigma(2)\varsigma(5)$$

Using an integer relation detection algorithm, in 1994 Bailey, Borwein and Girgensohn [2] showed experimentally that

(14) $$24\sum_{n=1}^{\infty}\frac{\left[H_n^{(1)}\right]^4}{(n+1)^2}=859\varsigma(6)+72\varsigma^2(3)$$

and this relationship was also conjectured by Coffey [5] in 2005.

Assuming that this is correct, Choi and Srivastava [3] deduced that

(15) $$480\sum_{n=1}^{\infty}\frac{H_n^{(1)}H_n^{(3)}}{n^2}=9281\varsigma(6)-30\varsigma^2(3)$$

(16) $$120\sum_{n=1}^{\infty}\frac{\left[H_n^{(1)}\right]^2 H_n^{(2)}}{n^2}=1741\varsigma(6)+330\varsigma^2(3)$$

The above identities are all well-known to the followers of the mysteries of the Euler sums: so what's new?

The only new thing to report is the simple observation that the coefficient of the zeta function with the largest argument is a **prime** number and the other coefficients are all **even** integers. This statement does not apply to (7) which is simply used to construct (9) and (10) which do indeed contain those features. The above identities contain a widely dispersed set of prime numbers $\{2,5,17,43,97,101,373,557,859,1741,2357,9281\}$. No doubt Mr. Darwin would consider that this was a selection of some sort!

We do however come across exceptions; for example in the same paper Choi and Srivastava [3] report that

(17) $$80\sum_{n=1}^{\infty}\frac{H_n^{(1)}H_n^{(3)}}{n^2}-20\sum_{n=1}^{\infty}\frac{\left[H_n^{(1)}\right]^4}{n^2}=731\varsigma(6)-75\varsigma^2(3)$$

where 731 is not a prime and 75 is odd (731 is however very close to the proximate prime of 733). They also report that

(18) $$24\sum_{n=1}^{\infty}\frac{\left[H_n^{(1)}\right]^4}{n^2}=979\varsigma(6)+72\varsigma^2(3)$$

where $979=11.89$ (again it is very close to the proximate prime of 977).

We may readily manufacture other prime number relationships as follows. First of all, subtracting (1) from (2) gives us

(19) $$\sum_{n=1}^{\infty}\frac{\left[H_n^{(1)}\right]^2}{n^2}-\sum_{n=1}^{\infty}\frac{H_n^{(1)}}{n^3}=3\varsigma(4)$$

We also find that (2)-2(1) gives us

(20) $$4\sum_{n=1}^{\infty}\frac{\left[H_n^{(1)}\right]^2}{n^2}-8\sum_{n=1}^{\infty}\frac{H_n^{(1)}}{n^3}=7\varsigma(4)$$

and 2(2)-(1) will result in

(21) $$8\sum_{n=1}^{\infty}\frac{\left[H_n^{(1)}\right]^2}{n^2}-4\sum_{n=1}^{\infty}\frac{H_n^{(1)}}{n^3}=29\varsigma(4)$$

Similarly, 2(6)+(3) produces

(22) $$24\sum_{n=1}^{\infty}\frac{\left[H_n^{(2)}\right]^2}{n^2}+24\sum_{n=1}^{\infty}\frac{\left[H_n^{(1)}\right]^2}{n^4}=11\varsigma(6)$$

and from 5(12)+6(17) we have

(23) $$240\sum_{n=1}^{\infty}\frac{H_n^{(1)}H_n^{(2)}}{n^3}+480\sum_{n=1}^{\infty}\frac{H_n^{(1)}H_n^{(3)}}{n^2}-120\sum_{n=1}^{\infty}\frac{\left[H_n^{(1)}\right]^4}{n^2}=3881\varsigma(6)$$

and 3881 is a prime number. We therefore have an additional "manufactured" set of primes $\{3,7,11,29,3881\}$.

This may of course be just an uncanny coincidence, but the author thinks that it may be otherwise. It would be useful if more data could be collected for higher order Euler sums of even weight (and, in this regard, the "generating functions" in equation (4.3.32) et seq. of [7] and [11, p.250] may be of assistance). Reference should also be made to Coffey's paper [4] and to the recent paper by Zheng [13].

There may just possibly be a connection with Euler's celebrated formula [12, p.1]

(24) $$\varsigma(s)=\sum_{n=1}^{\infty}\frac{1}{n^s}=\prod_{\rho}\left(1-\frac{1}{\rho^s}\right)^{-1}$$

where $\rho$ runs through all of the primes.

A short table of prime numbers is set out below for ease of reference (for those who cannot sleep, the first 10,000 are reported in [1]).

| | | | | | | | | | |
|---|---|---|---|---|---|---|---|---|---|
| 2 | 3 | 5 | 7 | 11 | 13 | 17 | 19 | 23 | 29 |
| 31 | 37 | 41 | 43 | 47 | 53 | 59 | 61 | 67 | 71 |
| 73 | 79 | 83 | 89 | 97 | 101 | 103 | 107 | 109 | 113 |
| 127 | 131 | 137 | 139 | 149 | 151 | 157 | 163 | 167 | 173 |
| 179 | 181 | 191 | 193 | 197 | 199 | 211 | 223 | 227 | 229 |
| 233 | 239 | 241 | 251 | 257 | 263 | 269 | 271 | 277 | 281 |
| 283 | 293 | 307 | 311 | 313 | 317 | 331 | 337 | 347 | 349 |
| 353 | 359 | 367 | 373 | 379 | 383 | 389 | 397 | 401 | 409 |
| 419 | 421 | 431 | 433 | 439 | 443 | 449 | 457 | 461 | 463 |
| 467 | 479 | 487 | 491 | 499 | 503 | 509 | 521 | 523 | 541 |
| 547 | 557 | 563 | 569 | 571 | 577 | 587 | 593 | 599 | 601 |
| 607 | 613 | 617 | 619 | 631 | 641 | 643 | 647 | 653 | 659 |
| 661 | 673 | 677 | 683 | 691 | 701 | 709 | 719 | 727 | 733 |
| 739 | 743 | 751 | 757 | 761 | 769 | 773 | 787 | 797 | 809 |
| 811 | 821 | 823 | 827 | 829 | 839 | 853 | 857 | 859 | 863 |
| 877 | 881 | 883 | 887 | 907 | 911 | 919 | 929 | 937 | 941 |
| 947 | 953 | 967 | 971 | 977 | 983 | 991 | 997 | 1009 | 1013 |
| 1019 | 1021 | 1031 | 1033 | 1039 | 1049 | 1051 | 1061 | 1063 | 1069 |
| 1087 | 1091 | 1093 | 1097 | 1103 | 1109 | 1117 | 1123 | 1129 | 1151 |
| 1153 | 1163 | 1171 | 1181 | 1187 | 1193 | 1201 | 1213 | 1217 | 1223 |
| 1229 | 1231 | 1237 | 1249 | 1259 | 1277 | 1279 | 1283 | 1289 | 1291 |
| 1297 | 1301 | 1303 | 1307 | 1319 | 1321 | 1327 | 1361 | 1367 | 1373 |
| 1381 | 1399 | 1409 | 1423 | 1427 | 1429 | 1433 | 1439 | 1447 | 1451 |
| 1453 | 1459 | 1471 | 1481 | 1483 | 1487 | 1489 | 1493 | 1499 | 1511 |
| 1523 | 1531 | 1543 | 1549 | 1553 | 1559 | 1567 | 1571 | 1579 | 1583 |
| 1597 | 1601 | 1607 | 1609 | 1613 | 1619 | 1621 | 1627 | 1637 | 1657 |
| 1663 | 1667 | 1669 | 1693 | 1697 | 1699 | 1709 | 1721 | 1723 | 1733 |
| 1741 | 1747 | 1753 | 1759 | 1777 | 1783 | 1787 | 1789 | 1801 | 1811 |
| 1823 | 1831 | 1847 | 1861 | 1867 | 1871 | 1873 | 1877 | 1879 | 1889 |
| 1901 | 1907 | 1913 | 1931 | 1933 | 1949 | 1951 | 1973 | 1979 | 1987 |
| 1993 | 1997 | 1999 | 2003 | 2011 | 2017 | 2027 | 2029 | 2039 | 2053 |

Donal F. Connon
Elmhurst
Dundle Road
Matfield
Kent TN12 7HD
dconnon@btopenworld.com